\documentclass[draft,reqno]{amsart}

\usepackage{amsmath}
\usepackage{amssymb}
\usepackage[all]{xy}
\usepackage{picinpar}
\usepackage{palatino}
\usepackage[usenames,dvipsnames,svgnames,table]{xcolor}
\usepackage{mathtools}

\def\eu{\mathfrak}
\def\ma{\mathbb}
\def\mc{\mathcal}

\def\p{{\eu p}_{\infty}}
\def\f{{\ma F}_q^{\ast}}
\def\F{{\ma F}_q}

\def\pK{\eu p}
\def\pL{\eu P}
\def\fin{\hfill\qed\bigskip}
\def\ge#1{#1_{\eu{gex}}}
\def\g#1{#1_{\eu{ge}}}
\def\*#1{#1^*}
\def\cicl#1#2{k(\Lambda_{{#1}^{#2}})}

\def\lra{\longrightarrow}
\def\poly{P_1^{\alpha_1}\cdots P_r^{\alpha_r}}
\def\polyn#1{P_1^{\alpha_{1,#1}}\cdots P_r^{\alpha_{r,#1}}}
\def\Ku#1#2{k\big(\sqrt[#1]{#2}\big)}

\newcommand{\ord}{\operatorname{ord}}
\newcommand{\Gal}{\operatorname{Gal}}
\newcommand{\lcm}{\operatornamewithlimits{lcm}}
\newcommand{\mcd}{\operatorname{gcd}}

\newcommand{\Hom}{\operatorname{Hom}}
\newcommand{\xbinom}{\genfrac(){.5pt}0}

\newcounter{bean}
\newcounter{2bean}

\def\l{
\begin{list}
{\rm{(\alph{bean}).-}}{\usecounter{bean}
\setlength{\labelwidth}{0.8in}
\setlength{\labelsep}{0.3cm}
\setlength{\leftmargin}{1cm}}}

\def\las{\begin{list}
	{{\rm {(\arabic{2bean})}}}{\usecounter{2bean}
\setlength{\labelwidth}{0.8in}
\setlength{\labelsep}{0.3cm}
\setlength{\leftmargin}{1cm}}}

\numberwithin{equation}{section}
\newtheorem{theorem}{Theorem}[section]
\newtheorem{proposition}[theorem]{Proposition}
\newtheorem{lemma}[theorem]{Lemma}

\newtheorem{remark}[theorem]{Remark}

\newtheorem{corollary}[theorem]{Corollary}

\title[Genus fields of Kummer extensions of rational function fields]
{Genus fields of Kummer extensions of rational function fields}
 
\author[M. Rzedowski]{Martha Rzedowski--Calder\'on}
\address{Departamento de Control Autom\'atico\\
Centro de Investigaci\'on y de Estudios Avanzados del I.P.N.}
\email{mrzedowski@ctrl.cinvestav.mx}

\author[G. Villa]{Gabriel Villa--Salvador}
\address{Departamento de Control Autom\'atico\\
Centro de Investigaci\'on y de Estudios Avanzados del I.P.N.}
\email{gvillasalvador@gmail.com, gvilla@ctrl.cinvestav.mx}



\date{May 13th., 2021}

\begin{document}

\begin{abstract}

In this paper we obtain the genus field of a general Kummer
extension of a global rational function field. 
We study first the case of a general Kummer extension
of degree a power of a prime. Then we prove that
the genus field of a composite of two abelian extensions
of a global rational function field
with relatively prime degrees is equal to the composite
of their respective genus fields. Our main result,
the genus of a general Kummer extension of a
global rational function field,
is a direct consequence of this fact.

\end{abstract}

\maketitle

\section{Introduction}\label{S1}

The theory of genus fields for number fields has been around
in several presentations for more than two hundred years. Nowadays
it may be used to study the ``easy'' part of the Hilbert Class
Field (HCF) of a finite extension of the field of rational numbers.
The genus field of a number field $K$ (over ${\ma Q}$) is a
subfield of $K_H$, the HCF of $K$. In this setting, $K_H$ is
canonically defined as the maximal unramified
abelian extension of $K$ so the definition of the genus field
$\g K$ is also canonically given.

We are interested in global function fields. In this case, there are
several possible definitions for the HCF of a global function field $K$,
depending on which aspect of $K$ we are interested in. The simple
generalization of defining the HCF as the maximal unramified
extension of $K$, has the inconvenient of being of infinite degree
due to extensions of constants. Since every prime in a global 
function field becomes eventually inert in the extensions of constants,
a convenient way to define the HCF of $K$ is to require full
decomposition of at least one prime. When $K_H$ is finally defined,
the genus field $\g K$ of $K$ over a subfield $K_0$ is defined as
the composite $KK'$ so that $KK'\subseteq K_H$ and $K'/K_0$ is
the maximal abelian extension satisfying this property. In other
words, $K'$ is the maximal abelian extension of $K_0$ contained
in $K_H$.

Our setting is the following. We fix a rational function field $k=\F(T)$
and we consider a finite Kummer extension $K$ of $k$. We set
$S$ as the set of the infinite primes of $K$ and let $K_H$
be the maximal unramified abelian extension of $K$ such that
the elements of $S$ decompose fully. Our main result is
Theorem \ref{T4.9}, where we give a complete explicit expression
of the genus field $\g K$ of $K$ over $k$.

The tools we will be using in this paper are the cyclotomic function
fields given by the Carlitz module. In this setting 
we found in \cite{MaRzVi2013}
a general expression of the
genus field of any field cyclotomic function
field $E$ using the ramification theory of Dirichlet characters 
given by Leopoldt in \cite{Leo53}. Once we had the general
cyclotomic case, we were able to give a general expression
for a finite abelian extension of $K$. In the same paper, using
these general descriptions, we gave explicit expressions for
Kummer cyclic extensions of prime degree and for Artin--Schreier
extensions of $k$. In the first case it was previously obtained by
Peng \cite{Pen2003} and in the second one by Hu and Li \cite{HuLi2010}.
We also gave the explicit expression of $\g K$ for cyclic
extensions of degree $p^n$, where $p$ is the characteristic
of $k$, using Witt vectors.

In \cite{BaMoReRzVi2018} we presented another description of
genus fields of finite
abelian extensions of $k$, which is much more transparent and general
than the one in \cite{MaRzVi2013}. In particular we gave general
explicit expressions of genus fields of arbitrary 
abelian $p$--extensions. The case of cyclic extensions
of degree $l^n$ with $l$ a prime number $l\neq p$ has not
been solved in general. When the extension of degree $l^n$
is a Kummer extension, the explicit description of $\g K$ 
is given in \cite{MoReVi2019}. For an
extension that is not a Kummer extension,
we still do not have an explicit description of $\g K$.

The main goal of this paper is to give the explicit expression
of $\g K$ when $K/k$ is a finite Kummer extension. The
main difficulty is that in general, for two finite extensions
$K_1,K_2$ of $k$, we have $\g{(K_1)}\g{(K_2)}\subseteq 
\g{(K_1K_2)}$, but the equality does not always hold.
We will prove that, when the degrees of $K_i/k$ are
relatively prime, we have equality. This is Theorem
\ref{P4.2}. The explicit expression of the genus field of
a general Kummer extension of $k$, is now a
direct consequence of this equality and the case of prime
power degree.

In our study of genus fields, the main obstructions are the
appearance of inertia in the composite of fields and in
the proper contention of genus fields mentioned above.
These are the main reasons why the $p$--case has been solved
but the $l$--case, with $l\neq p$ has not.
To conclude the finite abelian case, it remains to study the
general cyclic case of prime power degree, not necessarily
Kummer, and the composite of several of these cases.

This paper can be considered as the end of the problem first
studied by Peng \cite{Pen2003} where he found the genus
field of a cyclic Kummer extension of a rational global function
field of degree $l$. The next step was the case of a Kummer cyclic
extension of prime power degree. This problem was first studied
in \cite{BaRzVi2013} under the condition that the field was
contained in a cyclotomic function field and a strong restriction.
In \cite{MoReVi2019} the general Kummer cyclic extension
of prime power degree was settled. Finally, here we obtained
explicitly the genus field of an arbitrary finite
Kummer extension of a global rational function field. This
is Theorem \ref{T4.9}. This result is a consequence of 
Theorems \ref{T3.4} and \ref{T3.6} where the general 
Kummer extension of prime power degree is established,
and of Theorem \ref{P4.2}.

\section{Notations and general results}\label{S2}

For the general Carlitz--Hayes 
theory of cyclotomic function fields, we
refer to \cite[Ch. 12]{Vil2006} and \cite[Cap. 9]{RzeVil2017}.
For the results on genus fields of function fields we
refer to \cite{BaMoReRzVi2018,MaRzVi2013,MaRzVi2017}
 and \cite[Cap. 14]{RzeVil2017}.

We will be using the following notation.
Let $k=\F(T)$ be a global rational function field,
where $\F$ is the finite field of $q$ elements. Let
$R_T=\F[T]$ and let $R_T^+$ denote the set of the monic
irreducible elements of $R_T$. For $N\in R_T$, $\cicl N{}$ denotes
the $N$--th cyclotomic function field where $\Lambda_N$ is
the $N$--th torsion of the Carlitz module. For $D\in R_T$ we
denote $\*D:=(-1)^{\deg D} D$.

We will call a field $F$ a {\em cyclotomic function field} 
if there exists $N\in R_T$ such that $F\subseteq \cicl N{}$.

Let $N\in R_T$. The Dirichlet characters $\chi\bmod N$ are
the group homomorphisms $\chi\colon (R_T/\langle N\rangle)^*
\lra \*{\ma C}$. Given a group $X$ of Dirichlet characters
modulo $N$, the {\em field associated to $X$} is the fixed field $F=
\cicl N{}^H$ where $H=\bigcap_{\chi\in X}\ker \chi$. We say that $F$
corresponds to the group $X$ and that $X$ corresponds to $F$.
We have that $X\cong \Hom(\Gal(F/k),\*{\ma C})$. When $X$ is a
cyclic group generated by $\chi$, we have that the field associated
to $X$ is equal to $F=\cicl N{}^{\ker \chi}$ and we say that $F$
corresponds to $\chi$.

Given a cyclotomic function field $F$ with group of Dirichlet characters $X$,
the maximal cyclotomic extension of $F$ unramified at the finite prime
divisors is the field that corresponds to $Y:=\prod_{P\in R_T^+}
X_P$ where $X_P=\{\chi_P\mid \chi\in X\}$ and $\chi_P$ is the
$P$--th component of $\chi$, see \cite{MaRzVi2013}.
We have that the ramification index of $P\in R_T^+$ in $F/k$ equals
$|X_P|$.

When $X$ is a cyclic group generated by $\chi$, the maximal cyclotomic
extension unramified at the finite primes of $F$ is given by $M=F_1\cdots F_r$
where $\chi=\prod_{i=1}^r \chi_{P_i}$ and $F_i:=\cicl N{}^{\ker \chi_{P_i}}$.
The only finite ramified prime in each $F_i/k$ is $P_i$, $1\leq i\leq r$.

We denote the infinite prime of $k$ by $\p$. That is, $\p$ is the pole
divisor of $T$ and $1/T$ is a uniformizer for $\p$. We have that
the inertia group and the decomposition group of $\p$ in 
$\cicl N{}/k$ are both equal to $\f\subseteq (R_T/\langle N\rangle)^*$.
In particular the ramification index of $\p$ in $\cicl N{}/k$ is equal to
$q-1$ and we define the the {\em maximal real subfield of $\cicl N{}$}
by the fixed field $\cicl N{}^+:=\cicl N{}^{\f}$. We have that $\p$
decomposes fully in $\cicl N{}^+/k$ and the inertia degree of $\p$
in every cyclotomic function field is always equal to $1$.

We define the maximal real subfield of a cyclotomic function field $F$
as $F^+:=F\cap \cicl N{}^+$, where $F\subseteq \cicl N{}$. We have that
$\p$ decomposes fully in $F^+/k$ and $\p$ is totally ramified in $F/F^+$.

Given a finite extension $K/k$ and a finite non-empty set $S$
of prime divisors of $K$, the {\em Hilbert Class Field} $K_{H,S}$
of $K$ with respect to $S$ is defined as the maximal unramified extension
of $K$ such that every prime in $S$ decomposes fully in $K_{H,S}/K$.
We will always take $S$ as the set of prime divisors dividing $\p$
and we will denote $K_{H,S}$ simply as $K_H$. The {\em genus field}
$\g K$ of $K$ with respect to $k$ is the maximal extension of $K$ contained
on $K_H$ that is of the form $K k^*$, where $k^*/k$ is an abelian extension.
We will always choose $k^*$ as the maximal extension
with respect to this property.
In other words, $\g K$ is equal to $K k^*$, where is $k^*$ is the maximal
abelian extension of $k$ contained in $K_H$. In the particular case
when $K/k$ is abelian, $\g K$ is the maximal abelian extension of $k$
containing $K$ such that $\g K/K$ is unramified and each element
of $S$ decomposes fully in $\g K/K$.

Let $F$ be any cyclotomic function field. Then $\g F=M^+ F$, where
$M$ is the maximal cyclotomic extension of $F$ unramified at the finite
primes. That is, $M$ is the field associated to $Y=\prod_{P\in R_T^+}
X_P$, where $X$ is the group of Dirichlet characters associated to $F$
(see \cite{BaMoReRzVi2018}). We denote $M=\ge F$ and then
$\g F=\ge F^+ F$. We have that $\ge F/\g F$ is totally ramified
at every element of $S$ and if $e_{\infty}(\ge F|F)$ denotes the
ramification index of every element of $S$ in $\ge F/F$, then
$[\ge F:\g F]=e_{\infty}(\ge F|F)=e_{\infty}(\ge F|\g F)$. Therefore, to
obtain $\g F$ we need to compute a subextension of $\ge F$ of
degree $e_{\infty}(\ge F|F)$.

We will use both notations: $e_{P}(F|k)$ or $e_{F/k}(P)$ to denote
the ramification index of the prime $P$ of $k$ in $F$.

When $K/k$ is a finite abelian extension, it follows from
the Kronecker--Weber Theorem that 
there exist $N\in R_T$, $n\in{\ma N}\cup\{0\}$ and $m\in{\ma N}$
such that $K\subseteq {_n\cicl N{}_m}$, where for any $F$, $F_m:=F
{\ma F}_{q^m}(T)$ and $_nF=FL_n$, with $L_n$ 
the maximal subfield of $\cicl {1/T^n}{}$,
where $\p$ is totally and wildly ramified.
Then we define $E:=K {\mc M}\cap \cicl N{}$, where 
${\mc M}=L_nk_m$. Then
$\g K=\g E^H K$, where $H$ is the decomposition group of
the infinite primes in $K\g E/K$ (see \cite{BaMoReRzVi2018}).

One important result on ramification of tamely ramified extensions, is
the following.

\begin{theorem}[Abhyankar's Lemma]\label{T2.1}
Let $L/K$ be a separable extension of global function fields. Assume that
$L=K_1K_2$ with $K\subseteq K_i\subseteq L$, $1\leq i\leq 2$. Let
$\pK$ be a prime divisor of $K$ and $\pL$ a prime divisor in $L$
above $\pK$. Let $\pL_i:=\pL\cap K_i$, $i=1,2$. If at least one of the
extensions $K_i/K$ is tamely ramified at $\pK$, then
\[
e_{L/K}(\pL|\pK)=\lcm[e_{K_1/K}(\pL_1|\pK),e_{K_2/K}(\pL_2|\pK)],
\]
where $e_{L/K}(\pL|\pK)$ denotes the ramification index.
\end{theorem}

\begin{proof}
See \cite[Theorem 12.4.4]{Vil2006}.
\end{proof}

\begin{remark}\label{R2.2}{\rm{
Abhyankar's Lemma is valid for function fields, not
only for global function fields.
}}
\end{remark}

\section{Kummer extensions of prime power degree}\label{S3}

Let $l$ be a prime with $l^n$ dividing $q-1$ and let
$K/k$ be a Kummer extension of exponent $l^n$. Then
from Kummer theory we have that $K$ is the composite
$K=K_1\cdots K_s$ of linearly disjoint cyclic Kummer
extensions. More precisely, $K$ can be written as
\begin{gather*}
K=k\big(\sqrt[l^{n_1}]{\gamma_1 D_1},
\cdots, \sqrt[l^{n_s}]{\gamma_s D_s}\big)=
K_1\cdots K_s\\
\intertext{where}
K_{\varepsilon}=k\big(\sqrt[l^{n_{\varepsilon}}]{\gamma_{\varepsilon} 
D_{\varepsilon}}\big), \quad 1\leq \varepsilon \leq s.\\
\intertext{We also have}
\Gal(K/k)\cong \Gal(K_1/k)\times\cdots\times
\Gal(K_s/k)\cong C_{l^{n_1}}\times\cdots\times C_{l^{n_s}}
\end{gather*}
with $n=n_1\geq n_2\geq \cdots \geq n_s$,
$\gamma_{\varepsilon}\in \f$, $D_{\varepsilon}\in R_T$ monic,
and $K_{\varepsilon}=k(\sqrt[l^{n_{\varepsilon}}]
{\gamma_{\varepsilon} D_{\varepsilon}})$ a cyclic extension of $k$ of
degree $l^{n_{\varepsilon}}$ for $1\leq \varepsilon \leq s$.

Let $P_1,\ldots,P_r$ be the set of finite primes of $k$ ramified
in $K$ with $P_1,\ldots,P_r\in R_T^+$ distinct. 
Then we may assume that 
\[
D_{\varepsilon}=\polyn {\varepsilon}\quad\text{with}\quad 0\leq
\alpha_{j,\varepsilon}\leq l^{n_{\varepsilon}}-1,\quad 1\leq j\leq r,
\quad 1\leq \varepsilon\leq s.
\]
In fact, $\alpha_{j,\varepsilon}=0$ if and only if $P_{j}$ is 
unramified in $K_{\varepsilon}/k$.

Let $\alpha_{j,\varepsilon}=b_{j,\varepsilon} l^{a_{j,\varepsilon}}$
with $\gcd(b_{j,\varepsilon},l)=1$ when $\alpha_{j,
\varepsilon}\neq 0$ and let $\deg P_{j}=c_{j} l^{d_{j}}$ with
$\gcd(c_{j},l)=1$, $1\leq j\leq r$.

 For $x\in{\ma Z}$,
$v_l(x)$ denotes the valuation of $x$ at $l$. That is, $v_l(x)=\gamma$
if $l^{\gamma}|x$ and $l^{\gamma+1}\nmid x$. We write $v_l(0)=\infty$.

\subsection{The cyclotomic case}
First, we assume
that $K$ is contained in a cyclotomic function field,
more precisely in $\cicl {D_1\cdots D_s}{}$, and this is so
if and only if $\gamma_{\varepsilon}\equiv
(-1)^{\deg D_{\varepsilon}}\mod ({\*\F})^{l^{n_{\varepsilon}}}$ for $1\leq
\varepsilon\leq s$ (see \cite[Corolario 9.5.12]{RzeVil2017}). 
When $K_{\varepsilon}$ is contained in a
cyclotomic function field, we may assume that $K_{\varepsilon}=k(\sqrt[l^{n_{\varepsilon}}]{
D_{\varepsilon}^*})$, $1\leq \varepsilon\leq s$. 
Note that if $l^{n_{\varepsilon}}|\deg D_{\varepsilon}$ then 
$\Ku{l^{n_{\varepsilon}}}{\*{D_{\varepsilon}}}=\Ku{l^{n_{\varepsilon}}}{D_{\varepsilon}}$.

First consider $F=\Ku{l^n}{\*D}$, with $D=\poly$,
a Kummer cyclic extension of $k$. Let $X=
\langle\chi\rangle$ be the group of Dirichlet characters associated to
$F$. Note that for any $\nu\in {\ma N}$ 
relatively prime to $l$, the
field associated to $\chi^{\nu}$ is $F$ since $X=\langle\chi^{\nu}\rangle$.
This corresponds to the fact that $F=\Ku{l^{n}}{\*{(D^{\nu})}}$. 

When $D=P\in R_T^+$
we have that the character associated to $F$ is $\xbinom{}P_{l^{n}}$,
the Legendre symbol which is defined as follows: if $P$ is of degree
$d$, then for any $N\in R_T$ with $P\nmid N$, $N\bmod P\in
(R_T/\langle P\rangle)^*\cong {\ma F}_{q^d}^*$. Then 
$\xbinom NP_{l^{n}}$
is defined as the unique element of ${\ma F}_{q^d}^*$ such that
$N^{\frac{q^d-1}{l^{n}}}\equiv \xbinom NP_{l^{n}}\bmod P$. We have that 
$\xbinom {}P_{l^{n}}$ is the character associated to $\Ku {l^{n}}{\*P}$
(see \cite[Proposici\'on 9.5.16]{RzeVil2017}). Let us denote
$\chi_P=\xbinom {}P_{l^{n}}$. Then, for any $\nu\in
{\ma Z}$, $\chi_P^{\nu}$ is the character
associated to $\Ku {l^{n}}{(P^{\nu})^*}$.

Hence, if $\chi_D$ is the character associated
to $\Ku {l^n}{\*D}$, then $\chi_{D}=
\prod_{j=1}^r\chi_{P_{j}}^{\alpha_{j}}$.

In general, for a radical extension, we have:

\begin{theorem}\label{T3.3} Let $F=\Ku{m}{\gamma D}$ be a geometric separable
extension of $k$, with $\gamma\in\f$ and let $D=\poly\in R_T$. Then 
\begin{gather*}
e_{F/k}(P_{j})=\frac m{\gcd(\alpha_j,m)},\quad 1\leq j
\leq r\quad \text{and}\\
e_{\infty}(F|k):=e_{F/k}(\p)=\frac m{\gcd(\deg D,m)}.
\end{gather*}
\end{theorem}

\begin{proof} 
See \cite[\S 5.2]{MaRzVi2017}.
\end{proof}

As a consequence we obtain the following result for a cyclic
cyclotomic Kummer extension $F=k(\sqrt[l^n]{\*D})$. Let $X$
be the group of Dirichlet characters associated to $F$ and let
$Y=\prod_{P\in R_T^+} X_P$ be the group
associated to $M$, the maximal
cyclotomic extension of $F$ unramified at the finite primes. 
Let $P=P_{j}$, $X=X_P=\langle\chi_P\rangle$ and let $F_0$ be
the field associated to $X_P$.
Then $F_0$ is cyclotomic, $P$ is the only ramified finite prime in
$F_0/k$ and $P$ is tamely ramified in $F_0/k$. This implies that $F_0
\subseteq \cicl P{}$ and $\Gal(\cicl P{}/k)\cong C_{q^{d_P}-1}$ with
$d_P:=\deg P$. Therefore $F_0$ is the only field of degree $\ord(
\chi_P)=l^{\beta_P}$ over $k$. Since
$F_0/k$ is a Kummer extension, it follows that $F_0
=\Ku{l^{\beta_P}}{\*P}$. Then, we have:

\begin{theorem}\label{T3.1}
The maximal unramified at the finite primes
cyclotomic extension of $F=\Ku {l^n}{\*D}$ is
$M:=k(\sqrt[l^n]{(P_1^{\alpha_1})^*},\ldots, \sqrt[l^n]{(P_r^{\alpha_r})^*})$.
$\fin$
\end{theorem}

\begin{remark}\label{R3.2}{\rm{
Let $\alpha=l^a b$ with $\gcd(b,l)=1$ and $0\leq a < n$. Then $\Ku {l^n}{(P^{\alpha})^*}=
\Ku {l^{n-a}}{\*P}$. In particular, if $\alpha_{\varepsilon}
=l^{a_{\varepsilon}}b_{\varepsilon}$ with $\gcd(l,b_{\varepsilon})=1$,
$1\leq \varepsilon\leq r$, then 
\[
M=k(\sqrt[l^{n-a_1}]{P_1^*},\ldots, \sqrt[l^{n-a_r}]{P_r^*})=
F_1\cdots F_r,
\]
with $F_{\varepsilon}=\Ku {l^{n-a_{\varepsilon}}}{P_{
\varepsilon}^*}$, $1\leq \varepsilon\leq r$.
}}
\end{remark}

We give another proof of Theorem \ref{T3.1} using Abhyankar's Lemma. On the one
hand we have that 
\[
[M:k]=\prod_{P\in R_T^+}|X_P|=\prod_{j=1}^r
|X_{P_{j}}|= \prod_{j=1}^r e_{M/k}(P_{j})=\prod_{j=1}^r 
l^{n-a_{j}}.
\]
On the other
hand if $F_{j}=\Ku{l^{n-a_{j}}}{\*{(P_{j})}}$, by Abhyankar's Lemma, 
we have $F F_{j}/K$ is
unramified at every finite prime, so $F F_1\cdots F_r/F$ is unramified at
the finite primes and $F\subseteq F_1\cdots F_r$. Hence $F_1\cdots
F_r\subseteq \ge F$ and $[F_1\cdots F_r:k]=[M:k]$. Therefore $M=
F_1\cdots F_r$ as claimed.

Let $\alpha_{j}=b_{j} l^{a_{j}}$ with $\gcd(l,b_{j})
=1$ and $\deg P_{j}=l^{d_{j}}c_{j}$ with $\gcd(l,c_{j})=1$, 
$1\leq j\leq r$. Then
\begin{gather}
e_{P_j}(F|k)=l^{n-a_{j}},\nonumber\\
e_{\infty}(F|k)=\frac {l^n}{\gcd(l^n,\deg D)}:=l^t,\label{Ec3.1}\\
e_{\infty}(F_{j}|k)=\frac{l^{n-a_{j}}}
{\gcd(l^{n-a_{j}},\deg P_{j})}=
\frac{l^{n-a_{j}}}{l^{\min\{n-a_{j},d_{j}\}}}=
l^{n-a_{j}-\min\{n-a_{j},d_{j}\}}.\nonumber
\end{gather}

From Abhyankar's Lemma we obtain that
\[
e_{\infty}(M|k)=\lcm_{1\leq j\leq r}[l^{n-a_{j}-\min\{
n-a_{j},d_{j}\}}]:=l^m.
\]

Therefore $[M:\g F]=l^{m-t}$. To find $\g F$ we must find a subfield
$F\subseteq L$ of $M$ such that $L/F$ is unramified and $[M:L]=l^{m-t}$.
If $L$ is such field, we have $L\subseteq \g F$ since $\p$ is
unramified in $L/F$ and $L$, being cyclotomic, 
satisfies that $\p$ decomposes
fully in $L/F$. Since $[M:L]=[M:\g F]$ it follows that $L=\g F$.

Going back to a general cyclotomic Kummer extension $K$ of $k$,
where $K=k\big(\sqrt[l^{n_1}]{\*{D_1}},
\cdots, \sqrt[l^{n_s}]{\*{D_s}}\big)=K_1\cdots K_s$.
Let $P_1,\ldots,P_r$ be the ramified finite primes in $K/k$.
From Abhyankar's Lemma and (\ref{Ec3.1}), we have, for $P_j\in R_T^+$
\begin{gather}
e_{P_j}(K|k)=\lcm_{1\leq\varepsilon\leq s}
[e_{P_j}(K_{\varepsilon}|k)]=l^{\beta_j}\nonumber
\intertext{with}
\beta_j=\max_{1\leq \varepsilon\leq s}\{n_{\varepsilon}-
v_l(\alpha_{j,\varepsilon})\}=\max_{\substack{1\leq \varepsilon\leq s\\
b_{j,\varepsilon}\neq 0}} \{n_{\varepsilon}-a_{j,\varepsilon}\},\label{Ec3.5}
\intertext{and}
\begin{align}
l^t:&=e_{\infty}(K|k)=\lcm_{1\leq \varepsilon\leq s}
[e_{\infty}(K_{\varepsilon}|k)]=
\lcm_{1\leq \varepsilon\leq s}\Big[\frac{l^{n_{\varepsilon}}}{\gcd(
l^{n_{\varepsilon}},\deg D_{\varepsilon})}\Big]\nonumber\\
&=\lcm_{1\leq \varepsilon\leq s}\big[l^{n_{\varepsilon}-\min\{
n_{\varepsilon},v_l(\deg D_{\varepsilon})\}}\big],\label{Ec3.6}
\end{align}
\end{gather}
that is, $t=\max\limits_{1\leq\varepsilon\leq s}\big\{n_{\varepsilon}-
\min\{n_{\varepsilon},v_l(\deg D_{\varepsilon}\}\big\}$.

Let $X_{\varepsilon}=\langle\chi_{\varepsilon}\rangle$ be
the group of Dirichlet characters corresponding to $K_{\varepsilon}$,
$1\leq \varepsilon\leq s$. Let $\chi_{\varepsilon}=\prod_{P\in 
R_T^+} \chi_{\varepsilon,P}$ be the product of $\chi_{\epsilon}$ into
its $P$--components, $e_P(K_{\varepsilon}|k)=\ord(\chi_{\varepsilon,P})$.

In this way, we obtain that $X=\langle
\chi_1,\ldots,\chi_s\rangle$ is the group of Dirichlet characters
associated to $K$, $X=X_1\cdots X_s$. We have
\[
X_P=(X_1)_P\cdots (X_s)_P=\langle\chi_{1,P}\rangle\cdots
\langle\chi_{s,P}\rangle=\langle\chi_{\gamma_P,P}\rangle
\]
with $\ord(\chi_{\gamma_P,P})=\max_{1\leq\varepsilon\leq s}
\{\ord(\chi_{\varepsilon,P})\}=e_P(K|k)$ for $P\in R_T^+$.

Let $M$ be the field associated to $Y:=\prod_{P\in R_T^+}X_P$,
the maximal cyclotomic extension of $K$ unramified at every finite prime.
From (\ref{Ec3.5}) we obtain
\[
e_{P_{j}}(K|k)=l^{\beta_j}.
\]

Then $M=F_1\cdots F_r$ with $F_{j}=\Ku{l^{\beta_j}}{\*{P_{j}}}$
and from (\ref{Ec3.1}) we obtain
\begin{align}
l^m:&=e_{\infty}(M|k)=\max_{1\leq j\leq r}\{e_{\infty}(F_j|k)\}=
\max_{1\leq j\leq r}\Big\{\frac{l^{\beta_j}}
{\gcd(l^{\beta_j},\deg P_j)}\Big\}\nonumber\\
&=\max_{1\leq j\leq r}\big\{
l^{\beta_j-\min\{\beta_j,d_j\}}\big\},\label{Ec3.7}
\end{align}
so that $m=\max\limits_{1\leq j\leq r}\big\{\beta_j-\min\{\beta_j,d_j\}\big\}$.

The procedure to obtain $\g K$ is the following. We order $P_1,\ldots,P_r$
so that $n=\beta_1\geq \beta_2\geq \cdots\geq \beta_r
\geq 1$, that is, we fix $P_1,
\ldots, P_r$ in decreasing order of their ramification indexes in $K/k$. There is
at least one $F_i$ such that $e_{\infty}(F_i|k)=l^m$. We choose $i$ as the
largest index with this property. Then we will show that there exist some
powers $z_j$, $1\leq j\leq i-1$ such that 
$\p$ is unramified in $\Ku {l^{\beta_j}}{(P_jP_i^{z_j})^*}/k$.

For $j>i$ we have two cases, $Q_j=P_jP_i^{y_jl^{\varepsilon_j}}$ 
or $Q_j=P_j^{y_j}P_i^{l^{\varepsilon_j}}$ for some $y_j,\varepsilon_j
\in {\ma Z}$ such that if 
$F_j=\Ku {l^{\gamma_j}}{Q_j}$ for some $\gamma_j$,
satisfies that the ramification index of $F_j/k$ at $P_j$ is $l^{\beta_j}$,
at $P_i$ is less than or equal to $l^{\beta_i}$ and $\p$ is unramified.
The rest will follow taking the composite of all these fields and one
of the form $\Ku{l^{\xi_i}}{P_i^*}$ for some $\xi_i$.

The result for the cyclotomic prime power degree case is the following.

\begin{theorem}\label{T3.4} 
Let $K/k$ be a finite cyclotomic Kummer $l$--extension 
of $k$, $K=K_1\cdots K_s$,
$K_{\varepsilon}=\Ku{l^{n_{\varepsilon}}}{D_{\varepsilon}^*}$,
$D_{\varepsilon}\in R_T$ monic, $1\leq \varepsilon\leq s$ and
$\Gal(K/k)\cong C_{l^{n_1}}\times \cdots \times C_{l^{n_s}}$
with $n=n_1\geq n_2\geq \cdots \geq n_s$ and $l^n|q-1$.
Then $K\subseteq \cicl {D_1\cdots D_s}{}$.
Let $P_1,\ldots,P_r$ be the finite primes in $k$ ramified in
$K$ with $P_1,\ldots,P_r\in R_T^+$ distinct. Let 
\begin{gather*}
e_{P_j}(K|k)=l^{\beta_j},\quad
1\leq \beta_j\leq n, \quad 1\leq j\leq r,\quad\text{and}\quad
e_{\infty}(K|k)=l^t, \quad 0\leq t\leq n
\end{gather*}
given by {\rm{(\ref{Ec3.5})}} and {\rm{(\ref{Ec3.6})}} and let
$\deg P_j=c_jl^{d_j}$ with $\gcd(c_j,l)=1$, $1\leq j\leq r$.

We order $P_1,\ldots, P_r$ so that $n=\beta_1\geq
\beta_2\geq \ldots \geq \beta_r$.

The maximal cyclotomic extension $M$ of $K$, unramified 
at the finite primes, is
given by $M=\ge K=k\Big(\sqrt[l^{\beta_1}]{P_1^*},\ldots,
\sqrt[l^{\beta_r}]{P_r^*}\Big)$. Let $l^m=e_{\infty}(M|k)$
be given by {\rm{(\ref{Ec3.7})}}.

Choose $i$ such that $m=\beta_i-\min\{\beta_i,d_i\}$
and such that for $j>i$ we have $m>\beta_j-\min\{\beta_j,d_j\}$. That is, $i$ is the largest
index obtaining $l^m$ as the ramification index of $\p$.

In case $m=t$ we have $M=\g K=
\prod_{j=1}^r\Ku{l^{\beta_j}}{P_j^*}$. 

In case $m>t\geq 0$, we have $\min\{
\beta_i,d_i\}=d_i$ and $m=\beta_i-d_i$. Let $a,b\in{\ma Z}$ be such that
$a\deg P_i+b l^{n+d_i}=l^{d_i}=\gcd(l^{n+d_i},\deg P_i)$. Set $z_j=-a
\frac{\deg P_j}{l^{d_i}}=-ac_jl^{d_j-d_i}\in {\ma Z}$ for $1\leq j\leq i-1$.
For $j>i$, consider $y_j\in{\ma Z}$ with $y_j\equiv -c_i^{-1}c_j\bmod l^n
=-ac_j \bmod l^n$. Let
\[
E_j=
\begin{cases}
\Ku{l^{\beta_j}}{P_jP_i^{z_j}}&\text{if $j<i$},\\
\Ku{l^{d_i+t}}{\*{P_i}}&\text{if $j=i$},\\
\Ku{l^{\beta_j}}{P_jP_i^{y_jl^{d_j-d_i}}}&\text{if $j>i$ and $d_j\geq d_i$},\\
\Ku{l^{\beta_j+d_i-d_j}}{P_j^{l^{d_i-d_j}}{P_i^{y_j}}}&\text{if $j>i$ and $d_i> d_j$}.
\end{cases}
\]

Then $\g K=E_1\cdots E_{i-1}E_iE_{i+1}\cdots E_r$.
\end{theorem}

\begin{proof}
When $m=t$ it follows that $\g K=M=\prod_{j=1}^r F_j$, where
$F_j=\Ku{l^{\beta_j}}{\*{P_j}}$.

Assume $m>t\geq 0$. Then $d_i<\beta_i$ and $\beta_i-d_i=m$. 

For $j<i$ we have $\beta_j\geq \beta_i$ and
$\beta_j-d_j\leq \beta_j-\min\{\beta_j,d_j\}\leq m=\beta_i-d_i$. Therefore $d_j\geq
\beta_j-m= \beta_j-(\beta_i-d_i)=\beta_j-\beta_i+d_i\geq d_i$. 
In particular $d_i|\deg P_j$,
$1\leq j\leq i-1$. Since $\gcd(l^{n+d_i},\deg P_i)=l^{d_i}$,
there exist $a,b\in{\ma Z}$,
such that $a\deg P_i+bl^{n+d_i}=l^{d_i}$. Multiplying by 
$\deg P_j$ and dividing by $l^{d_i}$, $j<i$, we obtain 
\[
a\frac{\deg P_j}{l^{d_i}}\deg P_i+b\deg P_j l^{n}=\deg P_j.
\]
That is, $\deg P_j+z_j\deg P_i=b\deg P_j l^{n}$, where
$z_j= -a\frac{\deg P_j}{l^{d_i}}=-ac_j l^{d_j-d_i}$. 
Note that $\mcd(a,l)=1$.
We have that $l^n|\deg(P_jP_i^{z_j})$.
Let $E_j:=\Ku{l^{\beta_j}}{P_jP_i^{z_j}}$. It follows that 
\begin{gather*}
e_{\infty}(E_j|k)=1, \quad e_{P_j}(E_j|k)=l^{\beta_j}=
e_{P_j}(K|k)\quad \text{and}\\
e_{P_i}(E_j|k)=l^{\beta_j-v_l(z_j)}|l^{\beta_i}=e_{P_i}(K|k),
\end{gather*}
since $v_l(z_j)=d_j-d_i$ and $\beta_j-v_l(z_j)=
\beta_j-d_j+d_i\leq m+d_i=\beta_i$. In particular $E_j\subseteq \g K$.

Now consider $j>i$. Let $y_j\in{\ma Z}$ such that
$y_j\equiv -c_i^{-1}c_j\bmod l^n$. Since $ac_i\equiv 1\bmod
l^n$ we have $c_i^{-1}\equiv a \bmod l^n$.
This is possible since $\mcd(c_ic_j,l)=1$.
Note that $\mcd(y_j,l)=1$. 

First assume $d_j\geq d_i$. Let
\begin{gather*}
Q_j=P_jP_i^{y_j l^{d_j-d_i}}\quad \text{and}\quad E_j:=\Ku{l^{\beta_j}}{Q_j}.
\intertext{We have}
\deg Q_j=\deg P_j+y_j l^{d_j-d_i} \deg P_i=c_jl^{d_j}+y_jl^{d_j-d_i} c_i
l^{d_i}=l^{d_j}(c_j+y_jc_i).
\end{gather*}
It follows that $l^n|\deg Q_j$ and $e_{\infty}(E_j|k)=1$. 
On the other hand 
\[
e_{P_j}(E_j|k)=l^{\beta_j}=e_{P_j}(K|k)\quad \text{and}\quad
e_{P_i}(E_j|k)=l^{\beta_j-d_j+d_i-v_l(y_j)}.
\]
Since $\gcd(y_j,l)=1$, we have
$\beta_j-d_j+d_i\leq m+d_i=\beta_i$ and $e_{P_i}(E_j|k)|e_{P_i}(K|k)$.
It follows that $E_j\subseteq \g K$.

Now consider $d_i\geq d_j$. Set $Q_j=P_j^{l^{d_i-d_j}}P_i^{y_j}$
and let $E_j=\Ku{l^{\beta_j+d_i-d_j}}{Q_j}$. We have
\[
\deg Q_j=l^{d_i-d_j}\deg P_j+y_j\deg P_i=l^{d_i-d_j}c_jl^{d_j}+y_jc_il^{d_i}=
l^{d_i}(c_j+y_jc_i).
\]
Therefore $l^n|\deg Q_j$ and $e_{\infty}(E_j|k)=1$. Next,
we have 
\begin{gather*}
e_{P_j}(E_j|k)=l^{\beta_j+d_i-d_j-d_i+d_j}=l^{\beta_j}=e_{P_j}(E_j|k) 
\intertext{and}
e_{P_i}(E_j|k)=l^{\beta_j+d_i-d_j}|l^{\beta_i}=e_{P_i}(K|k).
\end{gather*}
It follows that $E_j\subseteq \g K$.

Therefore $L:=E_1\cdots E_{i-1}E_{i+1}\cdots E_r\subseteq \g K$.
Note that from Abhyankar's Lemma, we have
\begin{gather*}
e_{\infty}(L|k)=1,\\\
e_{P_j}(L|k)=e_{P_j}(K|k)=l^{\beta_j}, j\neq i
\intertext{and}
e_{P_i}(L|k)|e_{P_i}(K|k)=l^{\beta_i}.
\end{gather*}
In fact, we can give a direct argument to show that actually $e_{P_i}(L|k)=
e_{P_i}(K|k)=l^{\beta_i}$ (see Remark \ref{R3.5}).

For any $1\leq j\leq r$ let $I_j$ denote the inertia group of $P_j$ in
$M'/k$ where $M'$ is any subfield of $M$ containing $E_j$. For any such 
$M'$ we have $|I_j|=l^{\beta_j}$.

Let ${\mc J}:=\{j>i\mid d_i>d_j\}$ and ${\mc I}:=\{1,2,\ldots, i-1,i+1,\ldots,
r\}\setminus {\mc J}$. In other words, if $j\in{\mc I}$ then $E_j=\Ku{l^{\beta_j}}{
P_jP_i^{x_j}}$ for some $x_j\in {\ma Z}$. 

Set ${\mc I}=\{m_1,\ldots,m_u\}$. Let $F:=E_{m_1}\cap E_{m_2}$. Since
we have that
$P_{m_1}$ is fully ramified in $E_{m_1}$ and is unramified in $E_{m_2}$,
it follows that $F=k$ and therefore 
$[E_{m_1}E_{m_2}:k]=[E_{m_1}:k][E_{m-2}:k]$.
\[
\xymatrix{
E_{m_1}\ar@{-}[r]\ar@{-}[d]&E_{m_1}E_{m_2}\ar@{-}[d]\\
k=E_{m_1}\cap E_{m_2}\ar@{-}[r]&E_{m_2}
}
\]
Furthermore, $\Gal(E_{m_1}E_{m_2}/k)\cong \Gal(E_{m_1}/k)\times
\Gal(E_{m_2}/k)\cong I_{m_1}\times I_{m_2}$ and $(E_{m_1}E_{m_2})^{
I_{m_1}I_{m_2}}=k$. By induction, we obtain for $1\leq v\leq u$:
\las
\item $\big(E_{m_1}\cdots E_{m_{v-1}}\big)\cap E_{m_v}=k$,
\item $\big[E_{m_1}\cdots E_{m_{v}}:k\big]=[E_{m_1}:k]\cdots [E_{m_v}:k]$,
\item $\big(E_{m_1}\cdots E_{m_{v}}\big)^{I_{m_1}\cdots I_{m_v}}=k$,
\item $I_{m_1}\cdots I_{m_v}\cong I_{m_1}\times\cdots\times I_{m_v}$.
\end{list}

That is, for any $\mu\in {\mc I}$ we have
that $\Big(\prod_{j\in{\mc I}\setminus \{\mu\}}E_j\Big)\cap 
E_{\mu}=k$ since for any non-trivial subfield
of $A:=\prod_{j\in{\mc I}\setminus \{\mu\}}E_j$ at least one $P_j$ with $j\in
{\mc I}\setminus \{\mu\}$ is ramified in this subfield and $P_j$ is unramified
in $E_{\mu}$. In particular we have
\begin{gather}\label{Ec3.2}
\big[\prod_{j\in{\mc I}}E_j:k\big]=\prod_{j\in{\mc I}}[E_j:k].
\end{gather}

We also have that
\begin{gather}\label{Ec3.3}
\Big(\prod_{j\in{\mc I}}E_j\Big)\bigcap \Big(\prod_{j\in{\mc J}}E_j\Big)=k
\end{gather}
since in any non-trivial subfield of $\prod_{j\in{\mc I}}E_j$ at least
one $P_{\mu}$ with $\mu\in {\mc I}$ is ramified and $P_{\mu}$ is unramified
in $\prod_{j\in{\mc J}}E_j$. In other words,
\[
[L:k]=\Big[\prod_{j\neq i}E_j:k\Big]=\Big[\prod_{j\in{\mc I}}E_j:k\Big]
\Big[\prod_{j\in{\mc J}}E_j:k\Big].
\]

To compute $\Big[\prod_{j\in{\mc J}}E_j:k\Big]$ we order ${\mc J}$ as
follows. Write ${\mc J}=\{j_1,\ldots,j_s\}$ with $d_i-d_{j_1}\leq d_i-
d_{j_2}\leq\cdots \leq d_i-d_{j_s}$. We have that $E_j^{I_j}=\Ku{l^{d_i-
d_j}}{P_j^{y_j}}=\Ku{l^{d_i-d_j}}{P_j}$.

First we consider $E_{j_1}E_{j_2}$. We have $E_{j_1}\cap E_{j_2}=
C_1$ where we denote $C_u:=\Ku{l^{d_i-d_{j_u}}}{P_i}$, 
$j_u\in{\mc J}$, $1\leq u\leq s$. In fact, if $\Lambda:=E_{j_1}\cap 
E_{j_2}$, then $P_{j_1}$ is not ramified in $E_{j_2}$ and $
P_{j_2}$ is not ramified in $E_{j_1}$, thus,
the only ramified prime in
$\Lambda/k$ is $P_i$. Furthermore, $C_1\subseteq E_{j_1}$
and $C_2\subseteq E_{j_2}$ and $C_1=C_1\cap C_2\subseteq
E_{j_1}\cap E_{j_2}$. Now, $P_{j_1}$ is fully ramified in $E_{j_1}/
C_1$. In particular, if $C_1\subsetneqq C'\subseteq F_{j_1}$, $P_{
j_1}$ is ramified in $C'/C_1$ and since $P_{j_1}$ is unramified in
($E_{j_1}\cap E_{j_2})/C_1$, it follows that $E_{j_1}\cap E_{j_2}=C_1$.

Consider the following diagram
\[
\xymatrix{
& E_{j_1}\ar@{-}[r]\ar@{-}[d]&E_{j_1}E_{j_2}\ar@{-}[d]\\
&C_1\ar@{-}[r]\ar@{-}[dl]&E_{j_2}\\ k
}
\]
We have $[E_{j_1}E_{j_2}:k]=[E_{j_1}:C_1][E_{j_2}:C_1][C_1:k]=
\frac{[E_{j_1}:k][E_{j_2}:k]}{[C_1:k]}$.

Now we consider $E_{j_1}E_{j_2}E_{j_3}$. With an argument similar
to the one in the previous case, we have
$E_{j_1}E_{j_2}\cap E_{j_3}=C_2$. We consider the following diagram.
\[
\xymatrix{
& E_{j_1}E_{j_2}\ar@{-}[r]\ar@{-}[d]&E_{j_1}E_{j_2}E_{j_3}\ar@{-}[d]\\
&C_2\ar@{-}[r]\ar@{-}[dl]&E_{j_3}\\ k
}
\]
Thus $[E_{j_1}E_{j_2}E_{j_3}:k]=[E_{j_1}E_{j_2}:C_2][E_{j_3}:C_2][C_2:k]
=[E_{j_1}E_{j_2}:k][E_{j_3}:C_2]=
\frac{[E_{j_1}:k][E_{j_2}:k][E_{j_3}:k]}{[C_1:k][C_2:k]}$.

By induction, we obtain for $1\leq v\leq s$
\las
\item $\big(E_{j_1}\cdots E_{j_{v-1}}\big)\cap E_{j_v}=C_{v-1}$, $v\geq 2$,
\item $\big[E_{j_1}\cdots E_{j_{v}}:k\big]=[E_{j_1}:C_1]\cdots 
[E_{j_v}:C_v][C_v:k]=\big(\prod_{\mu=1}^v l^{\beta_{j_{\mu}}}\big)
l^{d_i-d_{j_v}}$,
\item $\big(E_{j_1}\cdots E_{j_{v}}\big)^{I_{j_1}\cdots I_{j_v}}=k$,
\item $I_{j_1}\cdots I_{j_v}\cong I_{j_1}\times\cdots\times I_{j_v}$.
\end{list}
That is,
\begin{gather}
\Big[\prod_{j\in{\mc J}}E_j:k\Big]=[E_{j_1}\cdots E_{j_s}:k]=
\Big(\prod_{u=1}^s [E_{j_u}:C_u]\Big)[C_s:k]=
\big(\prod_{j\in{\mc J}}l^{\beta_j}\big)l^{d_i-d_s}.\label{Ec3.4}
\end{gather}

From (\ref{Ec3.2}), (\ref{Ec3.3}) and (\ref{Ec3.4}) we obtain
\begin{gather*}
[L:k]=\Big(\prod_{\substack{j=1\\ j\neq i}}^r l^{\beta_j}\Big) l^{d_i-d_s}
\quad\text{and}\quad L\cap E_i=\Ku{l^{d_i-d_{j_s}}}{P_i}=C_s,
\intertext{where $E_i:=\Ku{l^{d_i+t}}{\*{P_i}}$. Now, we have that}
e_{P_i}(E_i|k)=l^{d_i+t}|l^{\beta_i}=l^{d_i+m}=e_{P_i}(K|k)
\intertext{and}
e_{\infty}(F_i|k)=\frac{l^{d_i+t}}{\gcd(l^{d_i},l^{d_i+t})}=l^{d_i+t-d_i}
=l^t=e_{\infty}(K|k).
\end{gather*}

It follows that $LE_i\subseteq \g K$ and 
\begin{align*}
[LE_i:k]&=[L:C_s][E_i:C_s]
[C_s:k]=[L:C_s][E_i:k]\\
&=\frac{\prod_{j=1}^rl^{\beta_j}}{l^{\beta_i-(d_i+t)}}
=\frac{[M:k]}{l^{m-t}}=\frac{[M:k]}{[M:\g K]}=[\g K:k].
\end{align*}
Therefore $\g K=LE_i=E_1\cdots E_{i-1}E_iE_{i+1}\cdots E_r$. 
\end{proof}

\begin{remark}\label{R3.5}{\rm{
In the notation of Theorem \ref{T3.4}, we have that
in the case $m>t$, $K\subseteq \g K=E_1\cdots E_{i-1} E_i E_{i+1}\cdots
E_r$. Therefore there exists $j\neq i$ such that $e_{P_i}(E_j|k)=l^{\beta_i}$.
However one wonders why. Here we give a direct proof. For $j<i$
we have $E_j=\Ku{l^{\beta_j}}{P_jP_i^{z_j}}$ with $z_j=-acl^{d_j-d_i}$.
Now, $d_j\geq d_i$ and $\beta_j\geq \beta_i$ ($j<i$) and $v_l(z_j)=
d_j-d_i$ since $\gcd(ac_j,1)=1$. Then $e_{P_i}(E_j|k)=
l^{\beta_j-v_l(z_j)}=l^{\beta_j-d_j+d_i}$. So we require that
for some $j<i$ it holds $\beta_j-d_j+d_i=\beta_i$, equivalently,
$\beta_j-d_j=\beta_i-d_i$.

From the definition of the index $i$ we have $\beta_j-d_j\leq m=\beta_i-d_i$
and for $j>i$ we have $\beta_j-d_j<m$.

Assume that $\beta_j-d_j<m$ for all $j\neq i$. From (\ref{Ec3.5}) we
obtain that $\beta_j=\max_{1\leq \varepsilon\leq s}\{n_{\varepsilon}-
a_{j,\varepsilon}\}$ where for convenience we choose $a_{j,\varepsilon}=
n$ in case $\alpha_{j,\varepsilon}=b_{j,\varepsilon}l^{a_{j,\varepsilon}}=0$,
that is, when $b_{j,\varepsilon}=0$ because in this way $n_{\varepsilon}
-a_{j,\varepsilon}\leq 0$ and the maximum can not be obtained in $\varepsilon$
since $1\leq \beta_j\leq n$.

Let $1\leq \mu\leq s$ be such that $\beta_i=\max_{1\leq \varepsilon\leq s}
\{n_{\varepsilon}-a_{i,\varepsilon}\}=n_{\mu}-a_{i,\mu}$ so that $m=\beta_i-
d_i=n_{\mu}-a_{i,\mu}-d_i$ and 
\begin{gather}\label{Ec3.8}
a_{i,\mu}+d_i=n_{\mu}-m.
\end{gather}
Set $\deg D_{\mu}=c_0l^{d_0}$ with $\gcd (c_0,l)=1$, that is, $v_l(\deg D_{\mu})
=d_0$. Since $D_{\mu}=\prod_{j=1}^rP_j^{\alpha_{j,\mu}}$ we have
\begin{gather}\label{Ec3.9}
\deg D_{\mu}=\sum_{j=1}^r \alpha_{j,\mu}\deg P_j=\sum_{j=1}^r b_{j,\mu}
l^{a_{j,\mu}}c_jl^{d_j}=\sum_{j=1}^r b_{j,\mu}c_jl^{a_{j,\mu}+d_j}.
\end{gather}

Fix $j\neq i$ and let $\beta_j=\max_{1\leq\varepsilon\leq s}\{n_{\varepsilon}-
a_{j,\varepsilon}\}\geq n_{\mu}-a_{j,\mu}$. Therefore $a_{j,\mu}\geq n_{\mu}
-\beta_j$ and $a_{j,\mu}+d_j\geq n_{\mu}-\beta_j+d_j=n_{\mu}-(\beta_j-d_j)
>n_{\mu}-m$. Thus from (\ref{Ec3.8}) we obtain
\begin{gather}\label{Ec3.10}
a_{j,\mu}+d_j>n_{\mu}-m\quad \text{for}\quad j\neq i\quad\text{and}\quad
a_{i,\mu}+d_i=n_{\mu}-m.
\end{gather}

From (\ref{Ec3.9}) and (\ref{Ec3.10}) it follows that
\[
d_0=a_{i,\mu}+d_i=n_{\mu}-m.
\]
On the other hand we have
\begin{align*}
t&=\max_{1\leq \varepsilon\leq s}\{n_{\varepsilon}-\min\{n_{\varepsilon},
v_l(\deg D_{\varepsilon})\}\geq n_{\mu}-\min\{n_{\mu},d_0\}\\
&=n_{\mu}-\min\{n_{\mu},n_{\mu}-m\}=n_{\mu}-(n_{\mu}-m)=m,
\end{align*}
that is, $t\geq m$ contrary to our assumption: $t<m$. Therefore,
there exists $j<i$ such that $\beta_j-d_j=m$ and $e_{P_i}(E_j|k)=l^{\beta_i}
=e_{P_i}(K|k)$.
}}
\end{remark}

\subsection{The general case of prime power degree}\label{S3.2}

Now consider $K/k$ a Kummer extension of degree a power of the
prime number $l$.
If $\Gal(K/k)$ is of exponent $l^n$, we have $\Gal(K/k)\cong
C_{l^{n_1}}\times\cdots\times C_{l^{n_s}}$ with $n=n_1\geq\cdots
\geq n_s$ and $l^n|q-1$. We have that $K$ is of the form
$K=k\big(\sqrt[l^{n_1}]{\gamma_1 D_1},\ldots, 
\sqrt[l^{n_s}]{\gamma_s D_s}\big)$ with $D_{\varepsilon}\in
R_T$ monic, $\gamma_{\varepsilon}\in
\f$, $1\leq\varepsilon\leq s$ and $K_{\varepsilon}=\Ku{l^{n_{\varepsilon}}}
{\gamma_{\varepsilon} D_{\varepsilon}}$.

Let $E=Kk_{l^n}\cap \cicl {D_1\cdots D_s}{}=
k\big(\sqrt[l^{n_1}]{D_1^*},\ldots, \sqrt[l^{n_s}]{D_s^*}\big)$.
From \cite{BaMoReRzVi2018} we have that $\g K=\g E^H K$ where $H$
is the decomposition group of the infinite primes in
either $K\g E/\g K$ or $K\g E/K$.
We also have that $KE/KE^{H_1}$ is an extension of constants where
$H_1=H|_E$ is the decomposition group of the infinite
primes in $KE/K$ (see
\cite{BaMoReRzVi2018}).

We have that $EK=E\big(\sqrt[l^{n_1}]{\gamma_1 D_1},\ldots, 
\sqrt[l^{n_s}]{\gamma_s D_s}\big)=E\big(\sqrt[l^{n_1}]{\varepsilon_1},\ldots, 
\sqrt[l^{n_s}]{\varepsilon_s}\big)$ where $\varepsilon_{j}=(-1)^{\deg D_{j}}
\gamma_{j}$, $1\leq j\leq s$ since 
$\sqrt[l^{n_j}]{\varepsilon_j}=\frac{\sqrt[l^{n_j}]{\gamma_j D_j}}{\sqrt[l^{n_j}]
{D_j^*}}$ and $\sqrt[l^{n_j}]{D_j^*}\in E$. In particular, $EK/E$ is an extension
of constants and the inertia degree of $\p$ in $EK/k$ is
$f=l^v$ where ${\ma F}_{q^{l^v}}=\F\big(\sqrt[l^{n_1}]{\varepsilon_1},\ldots, 
\sqrt[l^{n_s}]{\varepsilon_s}\big)$ since $E$ being cyclotomic satisfies
that the inertia degree of $\p$ in $E/k$ is $1$.

It follows that
\[
|H|=\frac{[\F\big(\sqrt[l^{n_1}]{\varepsilon_1},\ldots, 
\sqrt[l^{n_s}]{\varepsilon_s}\big):\F]}{\deg_K \p}=\colon l^u.
\]

Now $H_1=H|_E\subseteq I_{\infty}(E|k)$ where $I_{\infty}(E|k)$
denotes the inertia group of $\p$ in $E/k$ and also denotes
the inertia group of $\p$ in $\g E/k$. In particular $\g E^+=
\g E^{I_{\infty}(E|k)}$. We denote by ${\mc H}_1$ the
inertia group of $\p$ in $\g E/k$. From Theorem \ref{T3.4}
we have that $\p$ is unramified in $L=E_1\cdots E_{i-1}E_{i+1}\cdots
E_r$ and totally ramified in $\g E/L\Ku{l^{d_i}}{P_i^*}$.
Therefore $\g E^+=L \Ku{l^{d_i}}{P_i^*}$ and
\[
I_{\infty}(\g E|k)=\Gal\big(\g E/L\Ku{l^{d_i}}{P_i^*}\big).
\]
The group ${\mc H}_1$ is the subgroup of order 
$l^u$ of $I_{\infty}(\g E|k)$.

It follows that $\g E^{{\mc H}_1}=L\Ku{l^{d_i+t-u}}{P_i^*}$ and therefore 
\[
\g K=\g E^{{\mc H}_1} K=E_1\cdots E_{i-1} \Ku{l^{d_i+t-u}}{P_i^*}
E_{i+1}\cdots E_r K.
\]

Therefore we have proved the main theorem on Kummer extensions
of prime power degree over $k$.

\begin{theorem}\label{T3.6}
Let $K=k\big(\sqrt[l^{n_1}]{\gamma_1 D_1},\ldots, 
\sqrt[l^{n_s}]{\gamma_s D_s}\big)$ 
be a Kummer extension of $k$ of prime power degree
with $D_{\varepsilon}\in
R_T$ monic, $\gamma_{\varepsilon}\in
\f$, $1\leq\varepsilon\leq s$ and $K_{\varepsilon}=\Ku{l^{n_{\varepsilon}}}
{\gamma_{\varepsilon} D_{\varepsilon}}$. Let $E=
k\big(\sqrt[l^{n_1}]{D_1^*},\ldots, 
\sqrt[l^{n_s}]{D_s^*}\big)$.

With the notations of Theorem {\rm{\ref{T3.4}}}, with $E$
in the place of $K$, we have that 
\[
\g K=L\Ku{l^{d_i+t-u}}{P_i^*}K=E_1\cdots E_{i-1}E_{i+1}\cdots E_r
\Ku{l^{d_i+t-s}}{P_i^*}K
\]
where $l^u=\frac{[\F(\sqrt[l^{n_1}]{\varepsilon_1},\ldots, 
\sqrt[l^{n_s}]{\varepsilon_s}):\F]}{\deg_K \p}$
and $\varepsilon_j =(-1)^{\deg D_j}\gamma_j$, $1\leq j\leq s$. $\fin$
\end{theorem}

\section{General Kummer extensions}\label{S4}

Let $K/k$ be a Kummer extension. Let $G:=\Gal(K/k)$
be the Galois group of $K/k$. If $S_1,\ldots, S_h$ are
the various Sylow subgroups of $G$, $G\cong S_1\times
\cdots\times S_h$ and each $S_j$ is a group of
prime power order, say $|S_j|=l_j^{m_j}$, $1\leq j\leq h$.

We may write $K=K_1\cdots K_h$ with $\Gal(K_j/k)
\cong S_j$, $1\leq j\leq h$.

From Theorem \ref{T3.6}, we know each $\g {(K_j)}$, $1\leq j\leq h$.
The knowledge of $\g K$ would follow if we had $\g K=
\prod_{j=1}^h\g{(K_j)}$. However in general for any
two fields $L_1$ and $L_2$, we only have that
$\g{(L_1)}\g{(L_2)}\subseteq \g{(L_1L_2)}$ 
(see \cite{MaRzVi2017} for an example where
$\g{(L_1)}\g{(L_2)}\subsetneqq \g{(L_1L_2)}$ ). 

Now, in our case we have $[K_i:k]=l_i^{m_i}$ and $\gcd([K_i:k],
[K_j:k])=1$ for all $i\neq j$. The knowledge of $\g K$ is an
immediate consequence of the following result.

\begin{theorem}\label{P4.2}
Let $L_i/k$, $i=1,2$ be two finite abelian extensions with
$\gcd([L_1:k],[L_2:k])=1$. Then $\g{(L_1)}\g{(L_2)}=\g{(L_1
L_2)}$.
\end{theorem}

To prove Theorem \ref{P4.2}, we first prove:

\begin{proposition}\label{P4.3}
Let $L/k$ be a finite abelian extension and let $l$ be any prime number.
Then $l|[L:k]$ if and only if $l|[\g L:k]$.
\end{proposition}

\begin{proof}
It is clear that if $l|[L:k]$ then $l|[\g L:k]$, because $L\subseteq \g L$.

Now let $l|[\g L:k]$. First assume that $L
\subseteq \cicl N{}$. Let $X$ be the group of Dirichlet characters
associated to $L$ and let $Y=\prod_{P\in R_T^+} X_P$.
Let $\ge L$ be the field associated to $Y$. We have
$\g L\subseteq \ge L$. In fact $\g L=\ge L^+L$. Then $l|[\ge L:k]$.
Therefore there exists $\chi\in Y$ of order $l$. Let $\chi
=\prod_{P\in R_T^+}\chi_P$ with $\chi_P\in X_P$ and $\chi
\neq 1$, $\chi^l=1$.

In particular there exists $P\in R_T^+$ with $\ord(\chi_P)=l$. Let
$\varphi\in X$ with
$\varphi_P=\chi_P$. Now if $l\nmid [L:k]$, $\gcd([L:k],l)=1$. Let
$m=\ord(\varphi)$. Then $l\nmid m$ and $\varphi^m=\prod_{
P\in R_T^+}\varphi_P^m=1$. Hence $\varphi_P^m=\chi_P^m=1$
and $\chi_P^l=1$. It follows that $\chi_P=1$. This contradicts
that $\ord\chi_P)=l$. It follows that $l|[L:k]$.

For the general case, with the notation as above, we have
$\g L=\g{E^H}L$. If $l|[\g L:k]$ then $l|
[\g{E^H}:k]$ or $l|[L:k]$. If $l|[\g{E^H}:k]$ then $l|[\g E:k]$ so,
from the cyclotomic case, we obtain that $l|[E:k]|[L:k]$. Hence
$l|[L:k]$.
\end{proof}

\begin{corollary}\label{C4.4} When $L$ is cyclotomic,
we have $l|[L:k]$ if and only if $l|[\ge L:k]$.
\end{corollary}

\begin{proof}
It follows immediately from the proof of Proposition \ref{P4.3}.
\end{proof}

\begin{lemma}\label{L4.5}
Let $L_i$, $i=1,2$ be two cyclotomic fields.
Then, if $\gcd ([L_1:k],[L_2:k])=1$, 
we have $(L_1L_2)^+=L_1^+L_2^+$.
\end{lemma}

\begin{proof}
Since $L_i^+\subseteq L_i$, $i=1,2$, it follows that $\gcd([L_1^+:k],
[L_2^+:k])=1$.

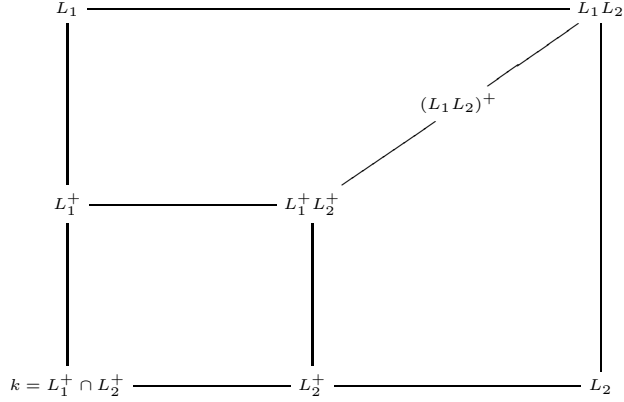
\begin{tiny}
\[
\xymatrix{
L_1\ar@{-}[dd]\ar@{-}[rrrr]&&&&L_1L_2\ar@{-}[dddd]\ar@{-}[dl]\\
&&&(L_1L_2)^+\ar@{-}[dl]\\
L_1^+\ar@{-}[rr]\ar@{-}[dd]&&L_1^+L_2^+\ar@{-}[dd]\\ \\
k=L_1^+\cap L_2^+\ar@{-}[rr]&&L_2^+\ar@{-}[rr]&&L_2
}
\]
\end{tiny}

We have that $e_i:=e_{\infty}(L_i|k)=[L_i:L_i^+]$, $i=1,2$.
Therefore $\gcd(e_1,e_2)=1$. It follows that $e_{\infty}(L_1L_2|
L_1^+L_2^+)=e_1e_2$ and $[L_1L_2:L_1^+L_2^+]=e_1e_2$.

That is, $L_1L_2/L_1^+L_2^+$ is totally ramified at $\p$.
On the other hand, $(L_1L_2)^+$ is the maximal subfield of
$L_1L_2$ where $\p$ decomposes fully in $(L_1L_2)^+/k$ and
$L_1L_2/(L_1L_2)^+$ is totally ramified at $\p$.
It follows that $e_{\infty}(
L_1L_2|(L_1L_2)^+)|e_1e_2=e_{\infty}(L_1L_2|L_1^+L_2^+)$. 
Thus $(L_1L_2)^+\subseteq L_1^+L_2^+\subseteq
(L_1L_2)^+$ and $(L_1L_2)^+=L_1^+L_2^+$. 
\end{proof}

\begin{remark}\label{R4.6}{\rm{
If $L_1$ and $L_2$ are cyclotomic it is possible that
$\gcd([L_1^+:k],[L_2^+:k])=1$ but $\gcd([L_1:k],
[L_2:k])\neq 1$.

For instance, let $P,Q\in R_T^+$ be distinct of degree $1$. Then
$[\cicl P{}:k]=q-1=[\cicl Q{}:k]$ and $\cicl P{}^+=\cicl Q{}^+=k$.
Thus, if $L_1=\cicl P{}$ and $L_2=\cicl Q{}$ we have 
$\gcd([L_1^+:k],[L_2^+:k])=\gcd(1,1)=1$ but $\gcd([L_1:k],
[L_2:k])=\gcd(q-1,q-1)=q-1>1$ for $q>2$.
}}
\end{remark}

In the cyclotomic case, the possibility $\g{(L_1)}\g{(L_2)}
\subsetneqq \g{(L_1L_2)}$ is related to the fact
that it may happen that $\ge{(L_1)}^+\ge{(L_2)}^+
\subsetneqq \ge{(L_1L_2)}^+$.
We will show that in our case we have equality.

\begin{corollary}\label{C4.7}
If $\gcd([L_1:k],[L_2:k])=1$ with $L_i\subseteq \cicl {N_i}{}$,
$i=1,2$, then $\g{(L_1L_2)}=\g{(L_1)}\g{(L_2)}$.
\end{corollary}

\begin{proof}
From Corollary \ref{C4.4} we obtain that 
$\gcd([\ge{(L_1)}:k],[\ge{(L_2)}:k])=1$. Now we have
$\g{(L_1)}=\ge{(L_1)}^+ L_1$ and $\g{(L_2)}=\ge{(L_2)}^+ L_2$.
Thus
\begin{align*}
\g{(L_1)}\g{(L_2)}&=\ge{(L_1)}^+ L_1\ge{(L_2)}^+ L_2=
\ge{(L_1)}^+\ge{(L_2)}^+ L_1L_2\\
&=(\ge{(L_1)}\ge{(L_2)})^+ (L_1L_2)=(\ge{(L_1L_2)})^+(L_1L_2)
=\g{(L_1L_2)}.
\end{align*}
\end{proof}

The proof of the following result is straightforward.

\begin{proposition}\label{P4.8}
Let $A,B$ and $C$ be global function fields such that
$B/A$ and $C/A$ are finite Galois extensions with
$\gcd([B:A],[C:A])=1$. Then $B\cap C=A$. Let $D=BC$.

If $\pK_A$ is a prime divisor of $A$ and $\pK_B,\pK_C,\pK_D$
satisfy $\pK_B\cap A=\pK_C\cap A=\pK_D\cap A=\pK_A$, then
\begin{gather*}
e_{\pK_A}(B|A)=e_{\pK_C}(D|C), \quad
f_{\pK_A}(B|A)=f_{\pK_C}(D|C), \quad
h_{\pK_A}(B|A)=h_{\pK_C}(D|C),\\
e_{\pK_A}(C|A)=e_{\pK_C}(D|B),\quad
f_{\pK_A}(C|A)=f_{\pK_C}(D|B),\quad
h_{\pK_A}(C|A)=h_{\pK_C}(D|B),
\end{gather*}
where $e,f$ and $h$ denote the ramification index, the inertia degree
and the decomposition degree respectively.
$\fin$
\end{proposition}

\subsubsection*{Proof of Theorem {\rm{\ref{P4.2}}}}

Let $L_i/k$ be two finite abelian extensions, $i=1,2$ such that
$\gcd([L_1:k],[L_2:k])=1$. Let $E_i=L_i{\mc M}\cap \cicl N{}$, $i=1,2$
where $L_i\subseteq {_n\cicl N{}_m}$, $i=1,2$. Let 
$L=L_1L_2$ and $E=E_1E_2$. Then $E=L{\mc M}\cap \cicl N{}$.

Now $\g{(L_i)}=L_i \g{(E_i)}^{H_i}$, $i=1,2$ with
$|H_i|=f_{\infty}(L_iE_i|L_i)$, $i=1,2$ and
$\g L=L \g E^H$ where $|H|=f_{\infty}(LE|L)$.
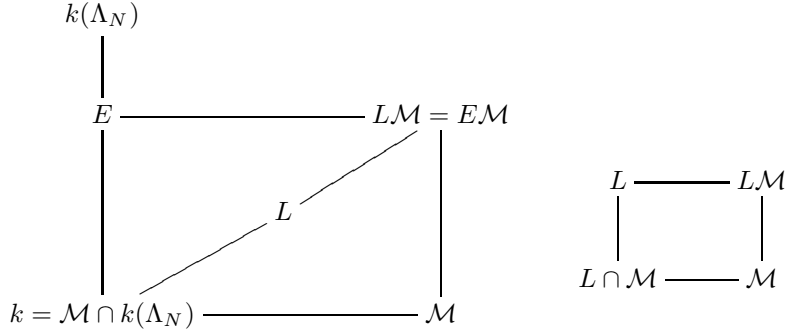
\[
\xymatrix{
\cicl N{}\ar@{-}[d]\\
E\ar@{-}[rr]\ar@{-}[dd]&&L{\mc M}=E{\mc M}\ar@{-}[dl]\ar@{-}[dd]\\
&L\ar@{-}[dl]\\
k={\mc M}\cap \cicl N{}\ar@{-}[rr]&&{\mc M}
}
\qquad
\xymatrix{
\\
\\
L\ar@{-}[r]\ar@{-}[d]&L{\mc M}\ar@{-}[d]\\
L\cap {\mc M}\ar@{-}[r]&{\mc M}
}
\]
We have $[E:k]=[L{\mc M}:{\mc M}]$ and $[L{\mc M}:{\mc M}]=[L:L\cap {\mc M}]|[L:k]$.
Therefore $[E:k]|[L:k]$. Analogously we have
$[E_i:k]|[L_i:k]$, $i=1,2$.

It follows that $\gcd([E_1:k],[E_2:k])=1$. From Proposition
\ref{P4.8} we obtain that
$|H|=|H_1||H_2|$. Then from Corollary \ref{C4.7} we obtain that
\begin{gather*}
[\g E:\g E^H]=|H|=|H_1||H_2|=[\g{(E_1)}:\g{(E_1)}^{H_1}]
[\g{(E_2)}:\g{(E_2)}^{H_2}]
\intertext{and we have, with $a=|H_1|$, $b=|H_2|$ and $ab=
|H_1||H_2|=|H|$, that}
\xymatrix{
\g{(E_1)}\ar@{-}[r]\ar@{-}[d]^a&\bullet\ar@{-}[d]^a
\ar@{-}[r]^{b\phantom{xxxxx}}&\g{(E_1)}\g{(E_2)}
=\g{(E_1E_2)}\ar@{-}[d]^a\ar@{-}[dl]_{ab}\\
\g{(E_1)}^{H_1}\ar@{-}[r]\ar@{-}[d] &\g{(E_1)}^{H_1}
\g{(E_2)}^{H_2}\ar@{-}[d]\ar@{-}[r]^b&\bullet
\ar@{-}[d] &\mcd(a,b)=1\\
k\ar@{-}[r]&\g{(E_2)}^{H_2}\ar@{-}[r]^b&\g{(E_2)}
}
\intertext{It follows that}
[\g{(E_1E_2)}:\g{(E_1)}^{H_1}\g{(E_2)}^{H_2}]=
[\g{(E_1)}:\g{(E_1)}^{H_1}][\g{(E_2)}:\g{(E_2)}^{H_2}]
=|H_1||H_2|=|H|
\end{gather*}
since $\g{(E_1E_2)}=\g{(E_1)}\g{(E_2)}$. On the other
hand, we have that
$[\g{(E_1E_2)}:\g{(E_1E_2)}^{H}]=|H|$. Since
$\g{(E_1)}^{H_1}\g{(E_2)}^{H_2}\subseteq \g{(E)}^{H}$,
we obtain that $\g{(E_1)}^{H_1}\g{(E_2)}^{H_2}=
\g{(E_1E_2)}^{H}$.
Therefore
\begin{align*}
\g{(L_1L_2)}&=\g{(E_1E_2)}^{H}(L_1L_2)=\g{(E_1)}^{H_1}
\g{(E_2)}^{H_2}L_1L_2\\
&=\g{(E_1)}^{H_1}L_1\g{(E_2)}^{H_2}
L_2=\g{(L_1)}\g{(L_2)}.
\tag*{$\fin$}
\end{align*}

As a consequence, we obtain the main result of the paper.

\begin{theorem}\label{T4.9}
Let $K/k$ be a finite Kummer extension of order
$n=l_1^{m_1}\cdots l_s^{m_s}$ with $l_1,\ldots,l_s$
different primes. Let $K=K_1\cdots K_s$ with $[K_j:k]=
l_j^{m_j}$, $1\leq j\leq s$. Then 
\[
\g K=\prod_{j=1}^s\g{(K_j)}
\]
where
each $\g{(K_j)}$, $1\leq j\leq s$ is computed in
Theorems {\rm{\ref{T3.4}}} and {\rm{\ref{T3.6}}}
\end{theorem}

\begin{proof}
From Theorem \ref{P4.2} we have $\g K=\g{(K_1)}\cdots
\g{(K_s)}$. The explicit description of each $\g{(K_i)}$ is the
content of Theorems \ref{T3.4} and \ref{T3.6}.
\end{proof}

\end{document}